\documentclass[10pt]{article} 
\begin{document}
\begin{flushleft}
\Large { \bf The Generation of Compatible Jacobi Tensors via Gauge
Transformations and its Applications} 
\newtheorem{definition}{Definition}
\newtheorem{Proposition}{Proposition}
\newtheorem{Lemma}{Lemma}
\newtheorem{Theorem}{Theorem}
\end{flushleft}
\begin{flushleft}
{\it Arthur SERGHEYEV} \footnote{e-mail:
 arthurser@imath.kiev.ua, arthur@apmat.freenet.kiev.ua}
\end{flushleft}
\begin{flushleft}
{\it Institute of Mathematics of the National Academy of Sciences of
Ukraine, \\ 
3 Tereshchenkivs'ka Street, Kyiv 4, \ Ukraine} \\
\end{flushleft}
\begin{minipage}{12cm}
\footnotesize
We present a simple way of generating the infinite set of Jacobi
tensors, compatible with a given one, via the "gauge transformations"
of the functions on Jacobi manifold.  We consider also some
applications of this result to the construction of bi-Hamiltonian
systems on Jacobi manifolds. 
\end{minipage} 

\vspace*{3mm}
\section{Introduction}

It is well known that Poisson manifolds play important role in
Hamiltonian mechanics and integrability \cite{o}. The suitable analytic
framework for these developments is usually referred as the
Poisson calculus \cite{vais}. The Poisson structure on manifold $M$
endows the space of functions $C^{\infty} (M)$ by the structure of
Lie algebra (via Poisson bracket), which respects the
usual structure of commutative algebra on $C^{\infty} (M)$, i.e. the
Leibnitz rule for the Poisson bracket takes place. Having rejected
the Leibnitz rule, we may obtain more general structures of Lie algebra
on $C^{\infty} (M)$, which are referred as local Lie algebras
\cite{kir}. The manifolds, equipped with such structures, are usually
called Jacobi manifolds. Namely, following Lichnerowicz \cite{lic},
we may define Jacobi manifold as follows
\begin{definition} Let $M$ be 
the manifold, on which are
defined a skew-symmetric 2-tensor $P$ and a vector $a$, which satisfy
\begin{eqnarray}
\lbrack P,P \rbrack = 2 a \wedge P,  \label{jac}\\
\lbrack P,a \rbrack = 0 \label{lie}.
\end{eqnarray}
Then the triple $(M,P,a)$ is called {\it Jacobi manifold}. 
\end{definition}

Here $[\cdot,\cdot]$ stands for the Schouten bracket (vide, e.g.,
\cite{o,lic}) and $\wedge$ stands for the exterior product of differential
forms.

The condition (\ref{lie}) means simply that the Lie derivative of $P$
with respect to $a$ vanishes. 

Provided $(M,P,a)$ is Jacobi manifold, the space $C^{\infty} (M)$
becomes a Lie algebra (more precisely, a {\it local} Lie algebra
\cite{kir}) with the commutator \cite{lic}
\begin{equation} \label{jb1}
\{ f,g \} =i(P) (df \wedge dg) + f i(a) dg - g i(a) df .
\end{equation} 

Note that in local coordinates $x^{k}, k=1, \dots, n={\rm dim}\; M$
in each chart of $M$ (\ref{jb1}) may be rewritten as  
\begin{equation} \label{jb}
\{ f,g \} = P^{kl}(x) \partial f/\partial x^{k}  \partial g/\partial x^{l}
+ a^{k} (x) (f \partial g/\partial x^{k} - g \partial f/\partial
x^{k}). 
\end{equation}

Here $P^{kl} = -P^{lk}$ and the Einstein's summation convention is
adopted with respect to the repeated indices $i,j,k,l, \dots$, which
are assumed to run from $1$ to $n$ (unless otherwise stated).\looseness=-1

The bracket (\ref{jb1}) is skew-symmetric by construction and
satisfies Jacobi identity in virtue of (\ref{jac}) and (\ref{lie}).

The pair $(P,a)$, satisfying (\ref{jac}) and (\ref{lie}), is usually
called the Jacobi tensor (or Jacobi structure \cite{lic}). Two Jacobi
tensors $(P_1,a_1)$ and $(P_2,a_2)$ are said to be {\it compatible} 
if $(P,a)=(P_{1}+P_{2}, a_{1} + a_{2})$ also is Jacobi tensor.

\section{The main result}
The straightforward computation shows that the following assertion
holds true:
\begin{Theorem} Let $(P,a)$ be Jacobi tensor on $M$ and $\varphi$ be an
arbirtary function from the space $C^{\infty} (M)$. Then 
\begin{equation} \label{new}
(\tilde P,
\tilde a) \equiv (\exp(\varphi) P,  \exp(\varphi) a - \exp (\varphi)
i(P) d\varphi)
\end{equation}
 is Jacobi tensor, compatible with $(P,a)$, and the map 
\[
\Psi : f
\mapsto f \exp (-\varphi), \quad f \in C^{\infty}(M), \quad (P,a)
\mapsto (\tilde P, \tilde a) 
\]
 is an isomorhism of local Lie algebras.
\end{Theorem}

Here $i(P) d\varphi$ stands for the vector with the components 
$P^{ij}(x) \partial \varphi/\partial x^{j}$.

As we see, the above isomorphism is constructed by means of the analog of
the well known gauge transformation in physics, $ f \mapsto f \exp
(i\chi)$, with setting $\chi = i \varphi$, where $i \equiv \sqrt{-1}$. 


Thus, the Jacobi tensors $(P,a)$ and $(\tilde P, \tilde a)$ are
equivalent and differ only by the "gauge" $\varphi$.

However, the compatibility of these Jacobi tensors allows one to
proceed in almost standard manner and try to define bi-Hamiltonian
systems on the Jacobi manifold in the following evident way:

Let $H \in C^{\infty}(M)$ be the Hamiltonian. Let us consider the
corresponding equations of motion
\begin{equation} \label{em1}
  df/dt = \{H,f \}_{(P,a)} \quad \mbox{for any}\: f \in C^{\infty}(M),
\end{equation}
where we fix explicitly the Jacobi tensor, which participates in the
definition of Jacobi bracket.

As we know, the Jacobi tensor $(\tilde P, \tilde a)$ is compatible
with $(P,a)$. Let us fix some non-constant\footnote{If
$\varphi$ is constant, the structure $(\tilde P, \tilde a)$ 
coincides with $(P,a)$ up to multiplication by the constant and
therefore is useless from the viewpoint of integrability of the
system (\ref{em1}).} function $\varphi$ in
(\ref{new}). Then, if there exists the function $H_{1} \in
C^{\infty}(M)$ such that 
\begin{equation} \label{bih}
 \{H,f \}_{(P,a)} = \{H_{1},f \}_{(\tilde P, \tilde a)},
\quad \mbox{for any}\: f \in C^{\infty}(M),
\end{equation}
the system (\ref{em1}) is bi-Hamiltonian, and one may 
construct the set of its integrals $I_{k}$ (which will be in
involution with respect to {\it both} brackets by construction
\cite{o}) in the standard way from the equation 
\begin{equation} \label{int}
\{I_{k-1},f \}_{(P,a)} = \{I_{k},f \}_{(\tilde P, \tilde a)}, \quad
k=1,2, \dots \quad \mbox{for any}\: f \in C^{\infty}(M).
\end{equation}
Two first terms in the sequence of $I_k$'s are $I_{0} =H$ and $I_{1}
= H_{1}$. 

Thus, if for any $k=1,2, \dots, k_{0}$ (in principle, the case $k_0
=\infty$ is not excluded) 
there exists the solution $I_{k}$ of (\ref{int}), and among
$I_k$, $k=0, 1, \dots, k_{0}$ one may choose a sufficient
number\footnote{On each symplectic leaf $S$ of $M$ this number is
$r=(1/2)\, {\rm dim}\;S$. By analogy with the theory of Poisson
structures, we mean here under the symplectic leaf a submanifold of
$M$, on which the bracket $\{,\}_{(P,a)}$ is non-degenerate.}   
of functionally independent ones, the system (\ref{em1}) will be
completely integrable in Liouville's sense.

\section{Discussion} 
We would like to note that our approach is quite different from the
usual scheme of 
proving the complete integrability in bi-Hamiltonian formalism, where the
key problem was the construction of the Poisson structure, compatible
with initial one. The key point is that we replace Poisson structure on
manifold by more general Jacobi structure. Then Theorem 1 always
guarantees the existence of the infinite set of Jacobi structures,
compatible with initial one. Hence, we have a very simple sufficient
condition for the Hamiltonian system (\ref{em1}) to be bi-Hamiltonian
and completely integrable: if for some non-constant
function $\varphi$ in (\ref{new}) one is able to find the sufficient
number of functionally independent solutions of (\ref{int}) for
different $k$, the system (\ref{em1}) will be completely integrable
in Liouville's sense, as it was already mentioned above.

The problem of existence of solutions of (\ref{int}) is now
under study and will be the subject of our following publications.

\medskip

I would like to thank Dr. Roman G. Smirnov for fruitful discussions 
about compatible Jacobi tensors and related topics.

\end{document}